\documentclass[12pt]{amsart}
\usepackage{verbatim}
\usepackage{eucal,url,amssymb,stmaryrd,enumerate,amscd,verbatim}
\usepackage{amsmath}
\usepackage[pagebackref,colorlinks=true,linkcolor=blue,citecolor=blue]{hyperref}
\usepackage{amsfonts}
\usepackage{amsthm}
\usepackage[margin=1in]{geometry}

\allowdisplaybreaks
\numberwithin{equation}{section}
\newtheorem{thrm}{Theorem}[section]
\newtheorem{lemma}[thrm]{Lemma}

\newtheorem{conv}[thrm]{Convention}

 1

\def\gr{\nabla f}
\def\bi{\nabla}

\newcommand{\vol}{\, Vol_{\eta}}

\begin{document}

\begin{abstract}
A quaternionic contact (qc) heat equation and the corresponding qc
energy functional are introduced. It is shown that the qc energy
functional is monotone non-increasing along the qc heat equation
on a compact qc manifold provided certain positivity conditions
are satisfied.
\end{abstract}

\title[A heat equation on a quaternionic contact manifold]{A heat equation on a quaternionic contact manifold}
\date{\today }
\author{S. Ivanov}
\address[Stefan Ivanov]{University of Sofia, Faculty of Mathematics and
Informatics, blvd. James Bourchier 5, 1164, Sofia, Bulgaria}
\address{and Institute of Mathematics and Informatics, Bulgarian Academy of
Sciences} \email{ivanovsp@fmi.uni-sofia.bg}
\author{A. Petkov}
\address[Alexander Petkov]{University of Sofia, Faculty of Mathematics and
Informatics, blvd. James Bourchier 5, 1164, Sofia, Bulgaria}
\email{a\_petkov\_fmi@abv.bg}

\maketitle

\tableofcontents


\setcounter{tocdepth}{2}

\section{Introduction}

We introduce a quaternionic contact (qc) heat equation and the
corresponding qc energy functional. The purpose of this paper is
to show that the qc energy functional is monotone non-increasing
along the qc heat equation on a compact qc manifold provided
certain positivity conditions are satisfied. In dimensions at
least eleven the positivity condition coincides with the
Lichnerowicz-type positivity condition used in \cite{IPV1,IPV3} to
derive a sharp  lower  bound for the first eigenvalue of the
sub-Laplacian an a compact qc manifold. In dimension seven, in
addition, we need to assume the positivity of the introduced in
\cite{IPV2} P-function.

It is well known that the sphere at infinity of a non-compact
symmetric space of rank one carries a natural
Carnot-Carath\'eodory structure, see \cite{M,P}. A quaternionic
contact (qc) structure, \cite{Biq1}, appears naturally as the
conformal boundary at infinity of the quaternionic hyperbolic
space. {Following Biquard, a quaternionic contact structure
(\emph{qc structure})  on a real (4n+3)-dimensional manifold $M$
is
a codimension three distribution $H$ (\emph{the horizontal distribution}) locally given as the kernel of a $%
\mathbb{R}^3$-valued one-form $\eta=(\eta_1,\eta_2,\eta_3)$, such
that the three two-forms $d\eta_i|_H$ are the fundamental forms of
a quaternionic Hermitian structure on $H$.
In other words, a quaternionic contact (qc) manifold $(M, g, \mathbb{Q})$ is a $4n+3$%
-dimensional manifold $M$ with a codimension three distribution
$H$ equipped with an $Sp(n)Sp(1)$ structure. Explicitly, $H$ is
the kernel of a local 1-form $\eta=(\eta_1,\eta_2,\eta_3)$ with
values in $\mathbb{R}^3$ together with a compatible Riemannian
metric $g$ and a rank-three bundle $\mathbb{Q}$ consisting of
endomorphisms of $H$ locally generated by three almost complex
structures $I_1,I_2,I_3$ on $H$ satisfying the identities of the
imaginary unit quaternions.

On a qc manifold one can associate a linear connection with
torsion preserving the qc structure, see \cite{Biq1}, which is
called the Biquard connection. } One defines the horizontal
Ricci-type tensor with the trace of the curvture  of the Biquard
connection,  called the qc Ricci tensor. This is a symmetric
tensor \cite{Biq1} whose trace-free part is determined by the
torsion endomorphism of the Biquard connection \cite{IMV} while
the trace part is determined by the scalar curvature of the
qc-Ricci tensor, called the qc-scalar curvature.

Let $(M,g,\mathbb{Q})$ be a compact 
qc
manifold. We consider \emph{the qc heat equation}
\begin{equation}\label{QCh}
\frac{\partial}{\partial t}u=-\Delta u,
\end{equation}
where $u(x,t):M\times[0,+\infty)\rightarrow\mathbb{R}$ is smooth
function and $\Delta:\mathcal{F}(M)\rightarrow\mathcal{F}(M)$ is
the sub-Laplacian on $M$. From now on, $u$ will be a positive
solution of \eqref{QCh}. We introduce the functions
$\varphi{\overset{def}=}-\ln u$ and $F{\overset{def}=}u^{\alpha},$
where $\alpha\in{\mathbb{R}}, \alpha\neq 0, \frac{1}{2}.$
\emph{The energy functional} for \eqref{QCh} is defined by
\begin{equation}\label{enf}
\mathcal{F}(\varphi)=\int_M|\nabla\varphi|^2\mathrm{e}^{-\varphi}\,Vol_{\eta}.
\end{equation}

Our main result follows.

\begin{thrm}\label{nonincr}
Let $(M,g,\mathbb{Q})$ be a compact $4n+3$-dimensional
quaternionic contact manifold and the Lichnerowicz type condition
\eqref{condm-app} holds, $L(X,X)\geq 0$ for any $X\in \Gamma (H)$.

\begin{itemize}
\item[i)] If $n>1$ then the energy functional \eqref{enf} is
monotone non-increasing along the qc heat equation \eqref{QCh}.
\item[ii)] In the case $n=1$ suppose in addition that the
$P-$function of any $F^{\frac{1}{2\alpha}},$ corresponding to a
(positive) solution $u$ of \eqref{QCh} is non-negative. Then the
energy functional \eqref{enf} is monotone non-increasing along the
qc heat equation \eqref{QCh}.
\end{itemize}
\end{thrm}
The Lichnerowicz
type assumption  cf. %
\eqref{Tcompnts}, \eqref{sixtyfour},
\begin{multline}  \label{condm-app}
L(X,X)=Ric(X,X)+\frac{2(4n+5)}{2n+1}T^0(X,X)+
\frac{6(2n^2+5n-1)}{(n-1)(2n+1)} U(X,X)\\=  2(n+2)Sg(X,X)+
\frac{4n^2+14n+12}{2n+1}T^0(X,X)
+\frac{4(n+2)^2(2n-1)}{(n-1)(2n+1)}U(X,X)\geq k_0 g(X,X),
\end{multline}
(the third term in the left-hand side is dropped if $n=1$) yields
a sharp lower bound of the first eigenvalue of the sub-Laplacian
when $n>1$ \cite{IPV1} while for $n=1$ one needs additional
assumption expressed in terms of the positivity of the
$P$-function defined in \cite{IPV3} to achieve the validity of the
same lower bound \cite{IPV3}. The $P$-function of a smooth
function $f$ is defined with the help
of the Biquard connection, the qc-scalar curvature and {the $Sp(n)Sp(1)$%
-components} of the torsion tensor see \eqref{e:def P} below.


\begin{conv}
\label{conven} \hfill\break\vspace{-15pt}

\begin{enumerate}[ a)]

\item We shall use $X,Y,Z,U$ to denote horizontal vector fields, i.e. $%
X,Y,Z,U\in H$.

\item $\{e_1,\dots,e_{4n}\}$ denotes a local orthonormal basis of the
horizontal space $H$.


\item The triple $(i,j,k)$ denotes any cyclic permutation of
$(1,2,3)$.

\item $s$ will be any number from the set $\{1,2,3\}$,
$s\in\{1,2,3\} $.
\end{enumerate}
\end{conv}

\textbf{Acknowledgments}  The authors thank Dimiter Vasssilev for
stimulating conversations. The research is partially supported by
Contract DFNI I02/4/12.12.2014 and by the Contract 195/2016 with
the University of Sofia `St.Kl.Ohridski'.

\section{Quaternionic contact manifolds}

Quaternionic contact manifolds were introduced in \cite{Biq1}. We also refer
to \cite{IMV} and \cite{IV} for further results and background.

\subsection{Quaternionic contact structures and the Biquard connection}

\label{ss:Biq conn}

A quaternionic contact (qc) manifold $(M, g, \mathbb{Q})$ is a $4n+3$%
-dimensional manifold $M$ with a codimension three distribution $H$ equipped
with an $Sp(n)Sp(1)$ structure. Explicitly, $H$ is the kernel of a local
1-form $\eta=(\eta_1,\eta_2,\eta_3)$ with values in $\mathbb{R}^3$ together
with a compatible Riemannian metric $g$ and a rank-three bundle $\mathbb{Q}$
consisting of endomorphisms of $H$ locally generated by three almost complex
structures $I_1,I_2,I_3$ on $H$ satisfying the identities of the imaginary
unit quaternions. Thus, we have $I_1I_2=-I_2I_1=I_3, \quad
I_1I_2I_3=-id_{|_H}$ which are hermitian compatible with the metric $%
g(I_s.,I_s.)=g(.,.)$ and the following compatibility conditions
hold $2g(I_sX,Y)\ =\ d\eta_s(X,Y)$.

On a qc manifold of dimension $(4n+3)>7$ with a fixed metric $g$
on  $H$ there exists a canonical connection defined in
\cite{Biq1}. Biquard also showed that there is a unique connection $%
\nabla$ with torsion $T$ and a unique supplementary subspace $V$ to $H$ in $%
TM$, such that:
\begin{enumerate}[i)]
\item $\nabla$ preserves the splitting $H\oplus V$ and the
$Sp(n)Sp(1)$ structure on $H$, i.e., $\nabla g=0, \nabla\sigma
\in\Gamma( \mathbb{Q})$ for a section
$\sigma\in\Gamma(\mathbb{Q})$, and its torsion on $H$ is given by
$T(X,Y)=-[X,Y]_{|V}$;
\item for $\xi\in V$, the endomorphism $T(\xi,.)_{|H}$ of $H$ lies in $%
(sp(n)\oplus sp(1))^{\bot}\subset gl(4n)$; \item the connection on
$V$ is induced by the natural identification $\varphi
$ of $V$ with $\mathbb Q$, 
$\nabla\varphi=0$.
\end{enumerate}


When the dimension of $M$ is at least eleven \cite{Biq1} also
described the supplementary \emph{vertical distribution} $V$,
which is (locally) generated by the so called \emph{Reeb vector
fields} $\{\xi_1,\xi_2,\xi_3\}$ determined by
\begin{equation}  \label{bi1}
\begin{aligned}
\eta_s(\xi_k)=\delta_{sk}, \qquad (\xi_s\lrcorner d\eta_s)_{|H}=0,
\qquad (\xi_s\lrcorner d\eta_k)_{|H}=-(\xi_k\lrcorner
d\eta_s)_{|H},
\end{aligned}
\end{equation}
where $\lrcorner$ denotes the interior multiplication.

If the dimension of $M $ is seven Duchemin shows in \cite{D} that if we
assume, in addition, the existence of Reeb vector fields as in \eqref{bi1},
then the Biquard result holds. \emph{Henceforth, by a qc structure in
dimension $7$ we shall mean a qc structure satisfying \eqref{bi1}}. This
implies the existence of the connection with properties (i), (ii) and (iii)
above.

The fundamental 2-forms $\omega_s$ of the quaternionic structure  are
defined by
\begin{equation*}
2\omega_{s|H}\ =\ \, d\eta_{s|H},\qquad \xi\lrcorner\omega_s=0,\quad \xi\in
V.
\end{equation*}
The torsion restricted to $H$ has the form $
T(X,Y)=-[X,Y]_{|V}=2\sum_{s=1}^3\omega_s(X,Y)\xi_s.
$

\subsection{Invariant decompositions}

Any endomorphism $\Psi$ of $H$ can be decomposed with respect to the
quaternionic structure $(\mathbb{Q},g)$ uniquely into four $Sp(n)$-invariant
parts 
$\Psi=\Psi^{+++}+\Psi^{+--}+\Psi^{-+-}+\Psi^{--+},$ 
where $\Psi^{+++}$ commutes with all three $I_i$, $\Psi^{+--}$ commutes with
$I_1$ and anti-commutes with the others two, etc. The two $Sp(n)Sp(1)$%
-invariant components \index{$Sp(n)Sp(1)$-invariant
components!$\Psi_{[3]}$} \index{$Sp(n)Sp(1)$-invariant
components!$\Psi_{[-1]}$} are given by $
\Psi_{[3]}=\Psi^{+++},\quad
\Psi_{[-1]}=\Psi^{+--}+\Psi^{-+-}+\Psi^{--+}. $
 These are the projections on the eigenspaces of the Casimir
operator $
\Upsilon =\ I_1\otimes I_1\ +\ I_2\otimes I_2\ +\ I_3\otimes I_3,
$
corresponding, respectively, to the eigenvalues $3$ and $-1$, see \cite{CSal}%
. Note here that each of the three 2-forms $\omega_s$ belongs to the
[-1]-component, $\omega_s=\omega_{s[-1]}$ and constitute a basis of the Lie
algebra $sp(1)$.

If $n=1$ then the space of symmetric endomorphisms commuting with all $I_s$
is 1-dimensional, i.e., the [3]-component of any symmetric endomorphism $\Psi
$ on $H$ is proportional to the identity, $\Psi_{[3]}=-%
\frac{tr\Psi}{4}Id_{|H}$.

\subsection{The torsion tensor}

The torsion endomorphism $T_{\xi }=T(\xi ,\cdot ):H\rightarrow H,\quad \xi
\in V$ will be decomposed into its symmetric part $T_{\xi }^{0}$ and
skew-symmetric part $b_{\xi },T_{\xi }=T_{\xi }^{0}+b_{\xi }$. Biquard
showed in \cite{Biq1} that the torsion $T_{\xi }$ is completely trace-free, $%
tr\,T_{\xi }=tr\,T_{\xi }\circ I_{s}=0$, its symmetric part has the
properties $T_{\xi _{i}}^{0}I_{i}=-I_{i}T_{\xi _{i}}^{0}\quad I_{2}(T_{\xi
_{2}}^{0})^{+--}=I_{1}(T_{\xi _{1}}^{0})^{-+-},\quad I_{3}(T_{\xi
_{3}}^{0})^{-+-}=I_{2}(T_{\xi _{2}}^{0})^{--+},\quad I_{1}(T_{\xi
_{1}}^{0})^{--+}=I_{3}(T_{\xi _{3}}^{0})^{+--}$.
The skew-symmetric part can be represented as $b_{\xi _{i}}=I_{i}U$, where $U
$ is a traceless symmetric (1,1)-tensor on $H$ which commutes with $%
I_{1},I_{2},I_{3}$. Therefore we have $T_{\xi _{i}}=T_{\xi _{i}}^{0}+I_{i}U$%
. When $n=1$ the tensor $U$ vanishes identically, $U=0$, and the
torsion is a symmetric tensor, $T_{\xi }=T_{\xi }^{0}.$

The two $Sp(n)Sp(1)$-invariant trace-free symmetric 2-tensors on $H$
\begin{equation} \label{Tcompnts}
T^0(X,Y)= g((T_{\xi_1}^{0}I_1+T_{\xi_2}^{0}I_2+T_{ \xi_3}^{0}I_3)X,Y) \
\text{ and }\ U(X,Y) =g(uX,Y)
\end{equation}
were introduced in \cite{IMV} and enjoy the properties
\begin{equation}  \label{propt}
\begin{aligned} T^0(X,Y)+T^0(I_1X,I_1Y)+T^0(I_2X,I_2Y)+T^0(I_3X,I_3Y)=0, \\
U(X,Y)=U(I_1X,I_1Y)=U(I_2X,I_2Y)=U(I_3X,I_3Y). \end{aligned}
\end{equation}
From \cite[Proposition~2.3]{IV} we have
\begin{equation}  \label{need}
4T^0(\xi_s,I_sX,Y)=T^0(X,Y)-T^0(I_sX,I_sY),
\end{equation}
hence, taking into account \eqref{need} it follows
\begin{multline}  \label{need1}
T(\xi_s,I_sX,Y)=T^0(\xi_s,I_sX,Y)+g(I_suI_sX,Y) \\
=\frac14\Big[T^0(X,Y)-T^0(I_sX,I_sY)\Big]-U(X,Y).
\end{multline}
\subsection{Torsion and curvature}

Let $R=[\nabla,\nabla]-\nabla_{[\ ,\ ]}$ be the curvature tensor of $\nabla$
and the dimension is $4n+3$. We denote the curvature tensor of type (0,4)
and the torsion tensor of type (0,3) by the same letter, $%
R(A,B,C,D):=g(R(A,B)C,D),\quad T(A,B,C):=g(T(A,B),C)$, $A,B,C,D
\in \Gamma(TM)$. {The \emph{qc-Ricci tensor} $Ric$,
\emph{normalized qc-scalar curvature} $S$ of the Biquard
connection are defined, respectively, by the following formulas
(cf. Convention 1.3),
$Ric(A,B)=\sum_{b=1}^{4n}R(e_b,A,B,e_b),\quad
8n(n+2)S=\sum_{a,b=1}^{4n}R(e_b,e_a,e_a,e_b).$ The qc-Ricci tensor
and the normalized qc-scalar curvature are determined by the
torsion of the Biquard connection a follows \cite{IMV}
\begin{equation} \label{sixtyfour}
\begin{aligned} & Ric(X,Y) =(2n+2)T^0(X,Y)+(4n+10)U(X,Y)+2(n+2)Sg(X,Y),
\\ &
T(\xi_{i},\xi_{j})=-S\xi_{k}-[\xi_{i},\xi_{j}]_{|H}, \qquad S =
-g(T(\xi_1,\xi_2),\xi_3).
\end{aligned}
\end{equation}
Note that for $n=1$ the above formulas hold with $U=0$.

Any 3-Sasakian manifold has zero torsion endomorphism, and the
converse is true if in addition the qc-scalar  curvature  is a
positive constant \cite{IMV}.

\subsection{The Ricci identities}

We  use repeatedly the Ricci identities of order two and three,
see also \cite{IV}. Let $f$ be a smooth function on the qc
manifold $M$ with horizontal gradient $\nabla f$ defined by
$g(\nabla f,X)=df(X)$. The sub-Laplacian of $f$ is $ \triangle
f=-\sum_{a=1}^{4n}\nabla^2f(e_a,e_a)$. We have the following Ricci
identities (see e.g. \cite{IMV,IV2})
\begin{equation*}  
\begin{aligned} & \nabla^2f
(X,Y)-\nabla^2f(Y,X)=-2\sum_{s=1}^3\omega_s(X,Y)df(\xi_s), \\ &
\nabla^2f (X,\xi_s)-\nabla^2f(\xi_s,X)=T(\xi_s,X,\nabla f),\\ &
\nabla^3 f (X,Y,Z)-\nabla^3 f(Y,X,Z)=-R(X,Y,Z,\nabla f) -
2\sum_{s=1}^3 \omega_s(X,Y)\nabla^2f (\xi_s,Z).
\end{aligned}
\end{equation*}
We also need  the qc-Bochner formula \cite[(4.1)]{IPV1}
\begin{multline}  \label{bohS}
\frac12\triangle |\nabla f|^2=|\nabla^2f|^2-g\left (\nabla
(\triangle f), \nabla f \right )+2(n+2)S|\nabla
f|^2+2(n+2)T^0(\nabla f,\nabla f) \\ +2(2n+2)U(\nabla f,\nabla f)+
4\sum_{s=1}^3\nabla^2f(\xi_s,I_s\nabla f). 
\end{multline}

\subsection{The horizontal divergence theorem}
Let $(M, g,\mathbb{Q})$ be a qc manifold of dimension $4n+3\geq 7$. For a
fixed local 1-form $\eta$ and a fixed $s\in \{1,2,3\}$ the form
$Vol_{\eta}=\eta_1\wedge\eta_2\wedge\eta_3\wedge\omega_s^{2n}$
is a locally defined volume form. Note that $Vol_{\eta}$ is independent of $%
s $ as well as the local one forms $\eta_1,\eta_2,\eta_3 $. Hence, it is a
globally defined volume form. The (horizontal) divergence of a horizontal
vector field/one-form $\sigma\in\Lambda^1\, (H)$, defined by $\nabla^*
\sigma =-tr|_{H}\nabla\sigma= -\nabla \sigma(e_a,e_a) $ 
supplies the integration by parts formula, \cite{IMV}, see also \cite{Wei},
\begin{equation*}  \label{div}
\int_M (\nabla^*\sigma)\,\, Vol_{\eta}\ =\ 0.
\end{equation*}

\subsection{The $P-$form}
We recall the definition of the P-form from \cite{IPV3}. Let
$(M,g,\mathbb{Q})$ be a compact quaternionic contact manifold of
dimension $4n+3$ and $f$ a smooth function on $M$.

For a smooth function $f$ on $M$ the $P-$form $P\equiv P_f \equiv
P[f]$ on $M$ is defined  by \cite{IPV3}
\begin{equation}
\begin{aligned} P_f(X) =&\nabla ^{3}f(X,e_{b},e_{b})+\sum_{t=1}^{3}\nabla
^{3}f(I_{t}X,e_{b},I_{t}e_{b})-4nSdf(X)+4nT^{0}(X,\nabla f)\\ &\hskip2.8in
-\frac{8n(n-2)}{n-1}U(X,\nabla f), \quad \text{if $n>1$},\\ P_f(X) =&\nabla
^{3}f(X,e_{b},e_{b})+\sum_{t=1}^{3}\nabla
^{3}f(I_{t}X,e_{b},I_{t}e_{b})-4Sdf(X)+4T^{0}(X,\nabla f),\ \text{if $n=1$}.
\end{aligned}  \label{e:def P}
\end{equation}

The $C-$operator is the fourth-order differential operator
independent of $f$ defined by
\begin{equation*}
Cf =-\nabla^* P_f=(\nabla_{e_a} P_f)\,(e_a).
\end{equation*}
We say that the $P-$function of $f$ is non-negative if its
integral exists and is non-positive
\begin{equation}  \label{e:non-negative Paneitz}
\int_M f\cdot Cf \, Vol_{\eta}= -\int_M P_f(\nabla f)\, Vol_{\eta}\geq 0.
\end{equation}
If \eqref{e:non-negative Paneitz} holds for any smooth function of
compact support we say that the $C-$operator is non-negative. It
turns out that the $C$-operator is non-negative on any compact qc
manifold of dimension at least eleven \cite{IPV3}.



One of the key identities which relates the P-function and the qc
Bochner formula \eqref{bohS} on a compact manifolds is  the next
identity, (dropping the last term when $n=1$), \cite[(3.4)]{IPV3}
\begin{multline}  \label{e:gr4}
\int_M\sum_{s=1}^3\nabla^2f(\xi_s,I_s\nabla f)\, Vol_{\eta} \\
=\int_M\Big[-\frac{1}{4n}P_f(\nabla f)-\frac{1}{4n}(\triangle f)^2-S|\nabla
f|^2+\frac{(n+1)}{n-1}U(\nabla f,\nabla f)\Big]\, Vol_{\eta}.
\end{multline}

\section{The QC heat equation and its energy functional}
The next lemma is crucial for the proof of our main result.
\begin{lemma}\label{deref}
Let $(M,g,\mathbb{Q})$ be a compact $4n+3$-dimensional
quaternionic contact manifold. Then the next formula holds
\begin{multline}\label{derf}
\alpha^2\frac{d}{d
t}\mathcal{F}(\varphi)=\frac{4\alpha}{3(1-2\alpha)}\int_MF^{\frac{1}{\alpha}-2}(\Delta
F)^2\,Vol_{\eta}\\+\frac{48n\alpha^2-2(16n-3)\alpha-3}{12(2n+1)\alpha^2}\int_MF^{\frac{1}{\alpha}-4}|\nabla
F|^4\,Vol_{\eta}+\frac{4(3-4\alpha)\alpha^2}{(2n+1)(1-2\alpha)}\int_MP_{F^\frac{1}{2\alpha}}(\nabla
F^\frac{1}{2\alpha})\,Vol_{\eta}\\-\frac{2n(3-4\alpha)}{3(n+2)(1-2\alpha)}\int_MF^{\frac{1}{\alpha}-2}L(\nabla
F,\nabla
F)\,Vol_{\eta}-\frac{4n(3-4\alpha)}{3(2n+1)(1-2\alpha)}\int_MF^{\frac{1}{\alpha}-2}p(F)\,Vol_{\eta}.
\end{multline}
\end{lemma}
In the formula \eqref{derf}, $P_{F^\frac{1}{2\alpha}}(\nabla
F^\frac{1}{2\alpha})$ is the $P-$function defined in \eqref{e:def
P} of $F^\frac{1}{2\alpha}$, $L(\nabla F,\nabla F)$ is the
left-hand side of the Lichnerowicz' type assumption
\eqref{condm-app} with $X:=\nabla F$ and
$$p(F){\overset{def}=}|\nabla^2F|^2-\frac{1}{4n}(\Delta
F)^2-\frac{1}{4n}\sum_{s=1}^{3}[g(\nabla^2F,\omega_s)]^2$$ is a
non-negative function on $M.$

\subsection{Proof of Lemma~\ref{deref}} The next relation
between the sub-Laplacians of $u$ and $\varphi$ holds
\begin{equation}\label{relsub}
\Delta u=-\frac{\Delta \varphi+|\nabla\varphi|^2}{\mathrm{e}^\varphi},
\end{equation}
which follows easily by the definitions of $\Delta$ and $\varphi$.
We get the formula
\begin{equation}\label{dtphi}
\frac{\partial}{\partial t}\varphi=-\Delta\varphi-|\nabla\varphi|^2,
\end{equation}
as a simply consequence of the definition of $\varphi,$ \eqref{QCh} and \eqref{relsub}. Further, the next chain of equalities holds
\begin{multline}\label{deriv1}
\frac{d}{dt}\mathcal{F}(\varphi)=\frac{d}{dt}\int_M\Big(-\Delta\varphi-\frac{\partial}{\partial t}\varphi\Big)u\,Vol_{\eta}=-\frac{d}{dt}\int_M\Delta\varphi u\,Vol_{\eta}+\frac{d}{dt}\int_M\Big(\frac{\partial}{\partial t}u\Big)\,Vol_{\eta}\\=-\int_M\Big[\Big(\Delta\frac{\partial}{\partial t}\varphi\Big)u+\Delta\varphi\frac{\partial}{\partial t}u\Big]\,Vol_{\eta}=-\int_M\Big(\frac{\partial}{\partial t}\varphi-\Delta\varphi\Big)\Delta u\,Vol_{\eta}\\=\int_M\mathrm{e}^{-\varphi}\Big[-2(\Delta\varphi)^2-3\Delta\varphi|\nabla\varphi|^2-|\nabla\varphi|^4\Big]\,Vol_{\eta},
\end{multline}
where we used \eqref{dtphi} for the first equality, the definition
of $\varphi$ for the second one, \eqref{QCh} and the divergence
theorem for the third equality. Finally, we took into account the
self-adjointness of the sub-Laplacian to obtain the fourth
equality and \eqref{relsub}, \eqref{dtphi} for the last one.

We need the next two identities:
\begin{equation}\label{connect}
|\nabla\varphi|^2=\alpha^{-2}F^{-2}|\nabla F|^2,\qquad\Delta\varphi=-\alpha^{-1}\Big(F^{-2}|\nabla F|^2+F^{-1}\Delta F\Big),
\end{equation}
which, substituted into \eqref{deriv1}, give
\begin{multline}\label{deriv2}
\alpha^2\frac{d}{dt}\mathcal{F}(\varphi)=-2\int_MF^{\frac{1}{\alpha}-2}(\Delta F)^2\,Vol_{\eta}\\+(3-4\alpha)\alpha^{-1}\int_MF^{\frac{1}{\alpha}-3}\Delta F|\nabla F|^2\,Vol_{\eta}+(-1+3\alpha-2\alpha^2)\alpha^{-2}\int_MF^{\frac{1}{\alpha}-4}|\nabla F|^4\,Vol_{\eta}.
\end{multline}

Next, we consider the (horizontal) vector field
$F^{\frac{1}{\alpha}-2}|\nabla F|^2,$ in order to deal with the
term $\int_MF^{\frac{1}{\alpha}-3}\Delta F|\nabla
F|^2\,Vol_{\eta}$ in \eqref{deriv2}. We get by some standard
calculations, using the divergence formula,
\begin{multline}\label{deriv3}
0=-\int_M\nabla^{*}\Big(F^{\frac{1}{\alpha}-2}|\nabla F|^2\Big)\,Vol_{\eta}\\=\int_Mg\Big(\nabla(F^{\frac{1}{\alpha}-2}\Delta F),\nabla F\Big)\,Vol_{\eta}-\int_MF^{\frac{1}{\alpha}-2}\Delta F\nabla^{*}\nabla F\,Vol_{\eta}\\=\int_MF^{\frac{1}{\alpha}-2}g\Big(\nabla(\Delta F),\nabla F\Big)\,Vol_{\eta}+\Big(\frac{1}{\alpha}-2\Big)\int_MF^{\frac{1}{\alpha}-3}\Delta F|\nabla F|^2\,Vol_{\eta}-\int_MF^{\frac{1}{\alpha}-2}(\Delta F)^2\,Vol_{\eta}.
\end{multline}
Integrate  the qc-Bochner formula \eqref{bohS} over the compact
$M$ and use \eqref{deriv3} to get
\begin{multline}\label{deriv4}
\Big(\frac{1}{\alpha}-2\Big)\int_MF^{\frac{1}{\alpha}-3}\Delta
F|\nabla
F|^2\,Vol_{\eta}\\=\int_MF^{\frac{1}{\alpha}-2}\Big[-\frac{1}{2}\Delta|\nabla
F|^2-|\nabla^2F|^2-2(n+2)S|\nabla F|^2-2(n+2)T^0(\nabla F,\nabla
F)\\-2(2n+2)U(\nabla F,\nabla
F)-4\sum_{s=1}^3\nabla^2F(\xi_s,I_s\nabla F)+(\Delta
F)^2\Big]\,Vol_{\eta}.
\end{multline}
The next step is to find some suitable representations of the two
terms $\int_MF^{\frac{1}{\alpha}-2}\Delta|\nabla F|^2\,Vol_{\eta}$
and
$\int_MF^{\frac{1}{\alpha}-2}\sum_{s=1}^3\nabla^2F(\xi_s,I_s\nabla
F)\,Vol_{\eta}.$ To deal with the first, we consider the
(horizontal) vector field $F^{\frac{1}{\alpha}-2}\nabla|\nabla
F|^2.$ We obtain the next sequence of equalities, using the
divergence formula and some standard calculations:
\begin{multline}\label{firstt}
0=-\int_M\nabla^{*}\Big(F^{\frac{1}{\alpha}-2}\nabla|\nabla F|^2\Big)\,Vol_{\eta}\\=\Big(\frac{1}{\alpha}-2\Big)\int_MF^{\frac{1}{\alpha}-3}g\Big(\nabla F,\nabla|\nabla F|^2\Big)\,Vol_{\eta}-\int_MF^{\frac{1}{\alpha}-2}\Delta|\nabla F|^2\,Vol_{\eta}\\=\Big(\frac{1}{\alpha}-2\Big)\int_MF^{\frac{1}{\alpha}-3}|\nabla F|^2\Delta F\,Vol_{\eta}-\Big(\frac{1}{\alpha}-2\Big)\Big(\frac{1}{\alpha}-3\Big)\int_MF^{\frac{1}{\alpha}-4}|\nabla F|^4\,Vol_{\eta}\\-\int_MF^{\frac{1}{\alpha}-2}\Delta|\nabla F|^2\,Vol_{\eta}.
\end{multline}
To get the third equality in \eqref{firstt} we used the identity
\begin{multline*}
0=\int_M\nabla^{*}\Big(F^{\frac{1}{\alpha}-3}|\nabla F|^2\nabla F\Big)\,Vol_{\eta}=-\int_MF^{\frac{1}{\alpha}-3}|\nabla F|^2\Delta F\,Vol_{\eta}\\+\int_MF^{\frac{1}{\alpha}-3}g\Big(\nabla F,\nabla|\nabla F|^2\Big)\,Vol_{\eta}+\Big(\frac{1}{\alpha}-3\Big)\int_MF^{\frac{1}{\alpha}-4}|\nabla F|^4\,Vol_{\eta}
\end{multline*}
in order to take an appropriate representation of the term $\int_MF^{\frac{1}{\alpha}-3}g\Big(\nabla F,\nabla|\nabla F|^2\Big)\,Vol_{\eta}.$

To handle the term
$\int_MF^{\frac{1}{\alpha}-2}\sum_{s=1}^3\nabla^2F(\xi_s,I_s\nabla
F)\,Vol_{\eta}$ we  use the next formula \cite[(3.12)]{IPV1}
\begin{equation}\label{intform}
\int_M\sum_{s=1}^3\bi^2f(\xi_s,I_s\gr)\vol=-\int_M\Big[4n\sum_{s=1}^3(df(\xi_s))^2
+\sum_{s=1}^3T(\xi_s,I_s\gr,\gr)\Big] \vol.
\end{equation}
Set $f:=F^\frac{1}{2\alpha}$ into \eqref{intform} to get after
some calculations that
\begin{multline}\label{secondt}
\int_MF^{\frac{1}{\alpha}-2}\sum_{s=1}^3\nabla^2F(\xi_s,I_s\nabla F)\,Vol_{\eta}\\=-\int_MF^{\frac{1}{\alpha}-2}\Big[4n\sum_{s=1}^3\Big(dF(\xi_s)\Big)^2+\sum_{s=1}^3T(\xi_s,I_s\nabla F,\nabla F)\Big]\,Vol_{\eta}
\end{multline}
Now, we substitute \eqref{firstt}, \eqref{secondt} in
\eqref{deriv4} and use the properties of the torsion tensor
\eqref{propt}, \eqref{need1} to obtain the identity
\begin{multline}\label{deriv5}
\frac{3}{2}\Big(\frac{1}{\alpha}-2\Big)\int_MF^{\frac{1}{\alpha}-3}\Delta F|\nabla F|^2\vol=\frac{1}{2}\Big(\frac{1}{\alpha}-2\Big)\Big(\frac{1}{\alpha}-3\Big)\int_MF^{\frac{1}{\alpha}-4}|\nabla F|^4\vol\\-\int_MF^{\frac{1}{\alpha}-2}\Big[|\nabla^2F|^2+2(n+2)S|\nabla F|^2+2nT^0(\nabla F,\nabla F)+4(n+4)U(\nabla F,\nabla F)\\-16n\sum_{s=1}^3\Big(dF(\xi_s)\Big)^2-(\Delta F)^2\Big]\vol.
\end{multline}
Substitute the right-hand side of \eqref{intform} into
\eqref{e:gr4} one obtains for $f:=F^{\frac{1}{2\alpha}}$ the
formula
\begin{multline}\label{intform1}
-4n\int_MF^{\frac{1}{\alpha}-2}\sum_{s=1}^3\Big(dF(\xi_s)\Big)^2\vol\\=\int_M\Big[-\frac{\alpha^2}{n}P_{F^{\frac{1}{2\alpha}}}(\nabla F^{\frac{1}{2\alpha}})-\frac{1}{4n}F^{\frac{1}{\alpha}-2}(\Delta F)^2+\frac{1}{2n}\Big(\frac{1}{2\alpha}-1\Big)F^{\frac{1}{\alpha}-3}\Delta F|\nabla F|^2\\-\frac{1}{4n}\Big(\frac{1}{2\alpha}-1\Big)^2F^{\frac{1}{\alpha}-4}|\nabla F|^4-F^{\frac{1}{\alpha}-2}\Big(S|\nabla F|^2-T^0(\nabla F,\nabla F)+\frac{2(n-2)}{n-1}U(\nabla F,\nabla F)\Big)\Big]\vol.
\end{multline}

It follows from the inequalities \cite[(4.6), (4.7)]{IPV1} the next representation of the norm of the horizontal Hessian:
\begin{equation}\label{hesrep}
|\nabla^2F|^2=\frac{1}{4n}(\Delta F)^2+\frac{1}{4n}\sum_{s=1}^3[g(\nabla^2F,\omega_s)]^2 +p(F),
\end{equation}
where $p(F)$ is a non-negative function on $M.$

Now, a substitution of \eqref{intform1} and \eqref{hesrep} in
\eqref{deriv5} give the identity

\begin{multline}\label{deriv6}
\int_MF^{\frac{1}{\alpha}-3}\Delta F|\nabla F|^2\vol=\frac{8\alpha^3}{(3n+2)(1-2\alpha)}\int_MP_{F^\frac{1}{2\alpha}}(\nabla F^\frac{1}{2\alpha})\vol\\+\frac{2n+1-2(3n+1)\alpha}{2(3n+2)\alpha}\int_MF^{\frac{1}{\alpha}-4}|\nabla F|^4\vol+\frac{(3+4n)\alpha}{2(3n+2)(1-2\alpha)}\int_MF^{\frac{1}{\alpha}-2}(\Delta F)^2\vol\\-\frac{2n\alpha}{(3n+2)(1-2\alpha)}\int_MF^{\frac{1}{\alpha}-2}\Big[2nS|\nabla F|^2+2(n+2)T^0(\nabla F,\nabla F)+\frac{4n(n+1)}{n-1}U(\nabla F,\nabla F)\Big]\vol\\-\frac{2n\alpha}{(3n+2)(1-2\alpha)}\int_MF^{\frac{1}{\alpha}-2}\Big[\frac{1}{4n}\sum_{s=1}^3[g(\nabla^2F,\omega_s)]^2+p(F)\Big]\vol.
\end{multline}
Note that we have the representation
\begin{multline}\label{represtor}
2nS|\nabla F|^2+2(n+2)T^0(\nabla F,\nabla F)+\frac{4n(n+1)}{n-1}U(\nabla F,\nabla F)\\=-S|\nabla F|^2+T^0(\nabla F,\nabla F)-\frac{2(n-2)}{n-1}U(\nabla F,\nabla F)+\frac{2n+1}{2(n+2)}L(\nabla F,\nabla F).
\end{multline}
Moreover, we obtain from the formula \cite[(4.12)]{IPV2}
\begin{multline*}
\int_M\Big[-S|\nabla f|^2+T^0(\nabla f,\nabla f)-\frac{2(n-2)}{n-1}U(\nabla
f,\nabla f)\Big]\,\,Vol_{\eta}=\int_M\Big[\frac{1}{4n}P_f(\nabla f)+\frac{1%
}{4n}(\triangle f)^2 \\
-\frac{1}{4n}\sum_{s=1}^{3}[g(\nabla^2f,\omega_s)]^2\Big]\,\,Vol_{\eta}
\end{multline*}
with $f:=F^{\frac{1}{2\alpha}}$ the next identity:
\begin{multline}\label{reprtor}
\int_MF^{\frac{1}{\alpha}-2}\Big[-S|\nabla F|^2+T^0(\nabla F,\nabla F)-\frac{2(n-2)}{n-1}U(\nabla F,\nabla F)\Big]\vol\\=\int_M\Big\{\frac{1}{4n}\Big[F^{\frac{1}{\alpha}-2}(\Delta F)^2-2\Big(\frac{1}{2\alpha}-1\Big)F^{\frac{1}{\alpha}-3}\Delta F|\nabla F|^2+\Big(\frac{1}{2\alpha}-1\Big)^2F^{\frac{1}{\alpha}-4}|\nabla F|^4\Big]\\+\frac{\alpha^2}{n}P_{F^\frac{1}{2\alpha}}(\nabla F^\frac{1}{2\alpha})-\frac{1}{4n}F^{\frac{1}{\alpha}-2}\sum_{s=1}^3[g(\nabla^2F,\omega_s)]^2\Big\}\vol.
\end{multline}

Taking into account \eqref{represtor} and \eqref{reprtor} in
\eqref{deriv6}, we get after some simple calculations
\begin{multline}\label{lastrep}
\frac{3(2n+1)}{2}\int_MF^{\frac{1}{\alpha}-3}\Delta F|\nabla F|^2\vol=\frac{8n+3-6(4n+1)\alpha}{8\alpha}\int_MF^{\frac{1}{\alpha}-4}|\nabla F|^4\vol\\+\frac{(2n+1)\alpha}{1-2\alpha}\int_MF^{\frac{1}{\alpha}-2}(\Delta F)^2\vol+\frac{6\alpha^3}{1-2\alpha}\int_MP_{F^\frac{1}{2\alpha}}(\nabla F^{\frac{1}{2\alpha}})\vol-\frac{2n\alpha}{1-2\alpha}\int_MF^{\frac{1}{\alpha}-2}p(F)\vol\\-\frac{n(2n+1)\alpha}{(n+2)(1-2\alpha)}\int_MF^{\frac{1}{\alpha}-2}L(\nabla F,\nabla F)\vol,
\end{multline}
which is the needed representation of the term $\int_MF^{\frac{1}{\alpha}-3}\Delta F|\nabla F|^2\vol.$

Finally, we substitute \eqref{lastrep} into \eqref{deriv2} to
obtain \eqref{derf}. This ends the proof of Lemma~\ref{deref}.

\subsection{Proofs of Theorem~\ref{nonincr}}
The polynomial
$h_n(\alpha){\overset{def}=}48n\alpha^2-2(16n-3)\alpha-3$ that
appears in the right-hand side of \eqref{derf} is non-positive for
$\alpha\in[\frac{16n-3-\sqrt{256n^2+48n+9}}{48n},\frac{16n-3+\sqrt{256n^2+48n+9}}{48n}].$
If we choose $\alpha\in[\frac{16n-3-\sqrt{256n^2+48n+9}}{48n},0)$
and suppose that the conditions (i) and (ii) of
Theorem~\ref{nonincr} hold, it is easy to see that  any summand
in the right-hand side of \eqref{derf} is non-positive, which
proofs Theorem~\ref{nonincr}.

\end{document}